\documentclass [12pt,a4paper,reqno]{amsart}
 \input amssymb.sty
\input xy
\xyoption{all}

\def\({\left(}
\def\){\right)}

\def\tvrp{\widetilde \vrp}

\def\Gm{\Gamma}
\def\mbR{\mathbb R}
\newcommand{\etype}[1]{\renewcommand{\labelenumi}{(#1{enumi})}}
\def\vrp{\varphi}
\def\pipeGS{{\underset{\operatorname{\, gs }}{\mid}}}

\def\lmod{\mathrel  \pipeGS   \joinrel \joinrel =}
\def\lmodg{\mathrel  \pipeGS   \joinrel \joinrel =}

\def\bbR{\mathbb R}
\def\bbN{\mathbb N}

\def\tG{\mathcal G}
\def\tT{\mathcal T}

\def\cG{\mathcal G}
\def\iff{\ \Leftrightarrow \ }

\def\ealph{\etype{\alph}}

\def\pSkip{\vskip 1.5mm \noindent}

\def\al{\alpha}

\def\gm{\gamma}

\newtheorem{thm}{Theorem} [section]
\newtheorem*{thm*}{Theorem}

\newtheorem{lem}[thm]{Lemma}

\newtheorem{prop}[thm]{Proposition}

\newtheorem*{claim*} {Claim}
\newtheorem*{theorem4.5'} {Theorem 4.5$'$}
\newtheorem{acknowledgment*}[thm] {Acknowledgment}
\newtheorem{example}[thm]{Example}
\newtheorem{sexample}[thm]{Subexample}

 \newtheorem*{remark*}{Remark}
 \newtheorem{defn}[thm]{Definition}

\newtheorem*{notation*} {Notation}
\newtheorem*{comment*} {Comment}

 \renewcommand{\sectionmark}[1]{}

\newcommand{\Cov}{\operatorname{Cov}}

\newcommand{\lm}{\lambda}

\def\ep{\varepsilon}
\def\tv{\tilde{v}}

\begin{document}

\title[A Glimpse at Supertropical Valuation Theory]
{A Glimpse at Supertropical Valuation Theory}
\author[Z. Izhakian]{Zur Izhakian}
\address{Department of Mathematics, Bar-Ilan University, 52900 Ramat-Gan,
Israel}
\email{zzur@math.biu.ac.il}
\author[M. Knebusch]{Manfred Knebusch}
\address{Department of Mathematics,
NWF-I Mathematik, Universit\"at Regensburg, 93040 Regensburg,
Germany} \email{manfred.knebusch@mathematik.uni-regensburg.de}
\author[L. Rowen]{Louis Rowen}
 \address{Department of Mathematics,
 Bar-Ilan University, 52900 Ramat-Gan, Israel}
 \email{rowen@macs.biu.ac.il}

\thanks{The research of the first and third authors has been supported  by the
Israel Science Foundation (grant No.  448/09).}

\thanks{The research of the second author was supported in part by
 the Gelbart Institute at
Bar-Ilan University, the Minerva Foundation at Tel-Aviv
University, the Department of Mathematics   of Bar-Ilan
University, and the Emmy Noether Institute at Bar-Ilan
University.}

\thanks{The research of the first author has been  supported  by the
Oberwolfach Leibniz Fellows Programme (OWLF), Mathematisches
Forschungsinstitut Oberwolfach, Germany.}


\thanks{This article is based on a lecture of the second author at
the International Conference and Humboldt Kolleg ``Fundamental
Structures of Algebra", at Constanta Romania, April 14--18, 2010,
in honor of the 70th birthday of Professor Serban Basarab.}

%
%

\date{\today}





\maketitle

\vskip 3mm \begin{center} \textsl{{To Serban A. Basarab with admiration, \\
on occasion of his seventieth birthday}}.
\end{center}
\vskip 15mm

 \baselineskip 14pt

\numberwithin{equation}{section}

%
%

\begin{abstract}
We give a short tour through major parts of a recent long paper
\cite{IKR1} on supertropical valuation theory, leaving aside
nearly all  proofs (to be found in \cite{IKR1}). In this way we
hope to give easy access to ideas of a new branch of so called
``supertropical algebra".

\end{abstract}

\section{Introduction}

 We will be much concerned  with
semirings. Recall that a \textbf{semiring} $R$ is a set $R$
equipped with addition and multiplication  such that both $(R, +
)$ and $(R \setminus \{ 0 \}, \cdot \, )$ are monoids, i.e.,
semigroups with a unit  element, $0$ and $1$ respectively, such
that multiplication distributes over addition  in the usual way.
In the present paper we always assume that multiplication (and, of
course, addition) is commutative. A \textbf{semifield} is a
semiring such that $(R \setminus \{ 0 \}, \cdot \, )$ is a group.
We give  two examples of semifields.

\begin{example}\label{exm:1.1} If $F$ is field then the set $R:= \sum F^2$
consisting of all  sums of squares in $R$ is a subsemiring of $F$,
and in fact a subsemifield since, if $$q : = a_1^2 + \cdots +
a_n^2$$ with $a_i \in F$ and $q \neq 0$, then $$ \frac 1 q = \(
\frac{a_1}{q} \)^2 + \cdots + \(\frac{a_n}{q}\)^2.$$

\end{example}

For the second example of semifields we need some preparation.  We
call a semiring $M$ \textbf{bipotent} if for any $a,b \in M$ the
sum $a+b$ is either $a$ or $b$. In this case we have a total
ordering $\leq$ on the set $M$, defined by
$$ a \leq b \iff a +b = b,  $$
as is easily checked. Clearly $0$ is the smallest element of $M$.
The ordering  is compatible with multiplication.
$$ a \leq b \ \Rightarrow \ ac \leq bc,  $$
and also with addition
$$ a \leq b \ \Rightarrow \ a+ c \leq b+c,  $$
for any  $a,b,c \in M$. We may state that
$$a+ b = \max(a,b).$$
Notice that bipotent semirings are very far away from those
semirings where addition is cancellative, i.e., $a+c = b+c
\Rightarrow a = b.$

\begin{example}\label{exm:1.2} Let $\Gm$ be a (totally) ordered
abelian  group, in multiplicative notation.  We add to $\Gm$ a new
element $0$ and extend the ordering of $\Gm$ to $M : = \Gm \cup \{
0 \}$ by declaring $0 < \gm$ for all $\gm \in \Gm$. We define
addition and multiplication on the set $M$ as follows:
$$\begin{array}{rll}
  x + y & = & \max(x,y) \quad \text{for } x,y \in M, \\[1mm]
  0 \cdot y & = & y \cdot 0 \ = \  0 \quad \text{for } y \in M, \\[1mm]
  x \cdot y & = &  \text{the given product in } \Gm, \text{ if }  x,y \in \Gm. \\
\end{array}$$
Clearly, $M$ is a bipotent semifield.
\end{example}

It is an easy exercise to check that in this way we obtain all
bipotent semifields $M$ from ordered abelian groups $\Gm$ in a
unique way ($M = \Gm \cup \{ 0 \}$, $\Gm= M \setminus \{ 0 \}$, $
\dots$). In short, \emph{bipotent semifields are the same objects
as ordered abelian  groups.}
\begin{sexample}\label{sexm:1.3} Take $\Gm = (\mbR, +)$, the
additive group of the real numbers with the standard ordering.
Since we switched to an additive notation, we denote the zero
element of the associated bipotent semiring $M$ now by $-\infty$.
Thus $M = \mbR \cup \{ -\infty\}$. Addition and multiplication on
$M$ are given by
$$ x \oplus y := \max(x,y), \qquad x\odot y : = x+y.$$
\end{sexample}

We refer to this bipotent semifield $\bbR \cup \{ -\infty \}$ and
related structures (e.g. the subsemiring $\bbR_{\geq 0} \cup \{
-\infty\}$) as the ``max-plus setting". It is used in tropical
geometry (e.g. \cite{Gathmann:0601322}, \cite{IMS}). \{In some
papers (e.g. \cite{SpeyerSturmfels2004}) an equivalent ``min-plus
setting" is used\}

The present authors feel that the max-plus setting is rather weak
for the needs of tropical geometry, and thus are driven by the
idea to develop a ``supertropical algebra", which should serve
tropical geometry better. Here the \textbf{supertropical
semirings}, to be defined below, occupy  a central place. The
prefix ``super" alludes to the fact that they are a sort of cover
of bipotent semirings.

There exist already supertropic results on polynomials
(\cite{IzhakianRowen2007SuperTropical},
\cite{IzhakianRowen2009Resulatants}), matrices
(\cite{IzhakianRowen2008Matrices},
\cite{IzhakianRowen2009Equations},
\cite{IzhakianRowen2010MatrixIII}) and, based on the supertropical
matrix theory, first steps of a supertropical linear algebra
\cite{IKR2}. And now supertropical valuation theory \cite{IKR1},
to which we refer  here.

\section {Supertropical predomains with pregiven ghost map}

\begin{defn}\label{defn:2.1} Let $R $ be a semiring. A \textbf{valuation $v$ on
$R$} is a map $v: R \to M$ into a bipotent semifield  $M$ with
$v(0)=0$, $v(1)= 1$, and, for any $a , b \in M$,
$$
\begin{array}{rll}
 v(a b ) & =&  v(a) \cdot v(b), \\[2mm]
 v(a + b ) & \leq  & v(a) +  v(b) \quad [= \max(v(a),v(b))] .
\end{array}$$
\end{defn}
If $R$ is a ring, this definition can be found in \cite[\S3
No.1]{Bo}, up to our change from ordered abelian groups
(additively written in \cite{Bo}) to bipotent semirings. If $R$ is
a field, we meet the classical Krull valuations.

\begin{defn}\label{defn:2.2} We call a valuation $v$ on the
semiring $R$ \textbf{strict}, if
$$\forall a,b \in R :  \qquad v(a+b) = v(a) + v(b),$$
(i.e., $v$ is a semiring homomorphism).

We call $v$ \textbf{strong}, if
$$\forall a,b \in R :  \qquad v(a) \neq v(b) \ \Rightarrow \   v(a+b) = v(a) + v(b).$$
\end{defn}

As is well known (at least for $R$ a field), every valuation $v$
on a ring $R$ is strong, but no valuation on $R$ is strict. But if
$R$ is  a semiring which is not ring, $v$ may be very well strict.

\begin{example}\label{exmp:2.3} (For readers with experience in real algebra.)
If $R = \sum F^2$ with $F$ a formal real field (cf. Example
\ref{exm:1.1}) and $w$ is a valuation on $F$, then the restriction
$w | R$ is strict iff the valuation $w$ is ``real", i.e., $w$ has
a formally real residue class field. In this way the real
valuations $w$ on $F$ correspond uniquely to the strict valuation
$v$ on $R$, provided the group $\Gm : = M \setminus \{ 0 \}$ is
$2$-divisible; we obtain $w$ back from $v$ by the formula
$$ w(a) = v(a^2)^{\frac{1}{2}}.$$
\end{example}

Since  any ordered abelian group can be enlarged to a
$2$-divisible ordered abelian group (even to a divisible ordered
abelian group) in a unique way, it is essentially a question of
preference, whether we study real valuations on fields or strict
valuations on sub-semifields. With the second route we leave the
cadre of classical algebra but have the possibility of transit to
semirings which  cannot be embedded into rings. For example we can
study the image of the ``total strict valuation map"
$$ R \to \prod M_v, \qquad a \mapsto (v(a_1)),$$
with $v$ running through all strict valuations $v: R \to M_v$ on
$R$. We do not pursue this line here, but only point out  that a
``semiring--approach"  is reasonable even for Krull valuations on
fields.

\section{Supertropical semirings}

We now define \textbf{supertropical semirings}. Such semirings
have first been constructed in a special case in \cite{I}, and
then defined in general in \cite{IzhakianRowen2007SuperTropical},
\cite{IzhakianRowen2008Matrices},  \cite{IKR1}. We follow the
approach of \cite{IKR1}, which has the advantage of being short,
but we refer the reader to the other papers to understand more on
the intuition behind these semirings.

\begin{defn}\label{defn:3.1} A semiring $U$ is
\textbf{supertropical} if the following four axioms ST1-ST4 hold,
which we state together with comments and side  definitions.
\begin{enumerate}
    \item[ST1:] The element  $e:= 1+1$ is idempotent, i.e., $1+1 =
    1+1+1+1$. \\ Thus $eU$ is an ideal of $U$ and is by itself a
    semiring.
\pSkip

     \item[ST2:] The semiring $eU$ is bipotent.

     Then the elements of $$ \cG := \cG(U) := eU \setminus \{0\} $$
are called the \textbf{ghost elements} of $U$. (In some sense also
$0$ is considered as a ghost element.) The map
$$ \nu_U:U \to eU, \qquad x \mapsto ex$$
is called the \textbf{ghost map} of $U$. It associates to each $x
\in U$ its \textbf{ghost} $ex$. (If $x \in eU$, then $x$ is it own
ghost.)

\pSkip
     \item[ST3:] If $ex < ey$ then $x+y = y$.

     \{Recall that $eU$ is totaly ordered, due to Axiom ST2.\}

\pSkip
    \item[ST4:] If $ex = ey$ then $x+y = ex$.

With Axiom ST4 we meet a principal idea of supertropical  algebra:
While in the usual tropical geometry the semirings are
\textbf{idempotent}, i.e., $x+ x =x$ for each $x$ in the semiring,
here $x+x$ is the ghost of $x$.

If $U$ is a supertropical semiring we call the elements of $$ \tT
:= \tT(U):= U \setminus eU$$ \textbf{tangible}. We then have  a
partition
$$ U = \tT \ \dot \cup \ \tG \ \dot\cup \ \{ 0 \},$$
and we remark that $\tG + \tG \subset \tG$.
\end{enumerate}

In the present paper we require  for supertropical semirings  one
more axiom, namely
\begin{enumerate}
    \item[ST5:] $\tT \cdot \tT \subset \tT, \quad \tG \cdot \tG \subset
    \tG$.
\end{enumerate}
\end{defn}
By this assumption we exclude only supertropical semirings which
are rather pathological and seldom of interest. (They are
sometimes needed for categorical reasons.)

We add three remarks for a supertropical semiring $U$.
\begin{enumerate}
    \item $U$ is bipotent iff $\tT$ is empty, \pSkip

    \item $\forall x\in U: \ ex =0 \ \Rightarrow \ x=0. $

    This is a consequence of ST4. We have $ex =e 0$, hence $x +0
    =0$.
\pSkip

    \item If $x_1, \dots, x_n \in U$  and $ x_1 + \cdots + x_n =
    0$, then all $x_i =0$.

    Indeed, we have $$ e x_1 + \cdots + e x_n =
    0,$$ and $ex_i \geq 0$ for each $i$. Since $U$ is totally ordered, it follows that each $ex_i=
    0$, hence $x_i = 0$.
\end{enumerate}

The second remark  indicates a special role of the zero element of
$U$. Informally it may be considered as both tangible and ghost.

We mention that there exists a completely  explicit construction
which gives us all supertropical semirings (with ST5), cf.
\cite[Construction 3.16]{IKR1}.

A basic intuition about ghost elements is that they are ``noise"
perturbing the tangible elements. This can be formulated as
follows:

\begin{defn}\label{defn:3.2} Given $x,y \in U$ we say that $x$
\textbf{surpasses $y$ by ghost}, and write $x \lmodg y$, if there
exists some $z \in eU$ with $x = y +z$.

We call the relation $\lmodg$ the \textbf{ghost surpassing
relation}, or \textbf{GS-relation}, for short.
\end{defn}

We state two remarkable properties of the GS-relation.
\begin{enumerate}
    \item  $\lmod$ is a partial ordering of the set $U$, which is
    compatible with multiplication, i.e.,
    $ x \lmod y \quad \text{implies} \quad  xz \lmod yz$
    for any $z \in U$. (The remarkable thing here is that $\lmod$
    is antisymmetric.)

    \item If $x \in \tT \cup \{ 0 \}$, $y\in U$, then $x \lmod y$
    implies $x=y$.
\end{enumerate}
    Thus if an element of $U$ is perturbed by adding a ghost, the resulting element can never be tangible.

\section{Supervaluatons}
We now introduce \emph{supervaluations}.

\begin{defn}\label{defn:4.1} A \textbf{supertropical semifield} is
a supertropical semiring $U$ for which the monoids $(\tT(U), \cdot
\;)$ and $(\tG(U), \cdot \;)$ are groups.
\end{defn}

Here we have to apologize for an inconsistency of language: The
ghost elements of a supertropical semifield are not invertible in
$U$ but only in $\tG(U)$. Thus $U$ is not a semifield as defined
in \S1, only $M:=eU$ is a semifield.

\begin{defn}\label{defn:4.2} A \textbf{supervaluation} on a
semiring $R$ is a map $\vrp: R \to U$ from $R$ to a supertropical
semifield $U$ with   $\vrp(0)=0$, $\vrp(1)= 1$, and, for any $a ,
b \in R$,
$$
\begin{array}{rll}
 \vrp(a b ) & =&  \vrp(a) \cdot \vrp(b), \\[2mm]
 e\vrp(a + b ) & \leq  & e\vrp(a) +  e\vrp(b).
\end{array}$$
\end{defn}

If $\vrp: R \to U$ is a supervaluation, then the map
$$ v: R \to M = eU, \qquad v(a):= e \vrp(a)$$
clearly is a valuation (as defined in \S2). We say that $\vrp$
\textbf{covers} the valuation $v$, and write $v = e \vrp$.

Starting with a valuation $v: R \to U$ with values in some
bipotent semifield $M$ we usually have very many supervaluations
$\vrp: R \to U$ covering $v$, where $U$ runs through the class of
all supertropical semifields with $M \subset U$ and $eU = M$. We
obtain a hierarchy between these supervaluations by a relation of
``dominance", to be explained now.

\begin{lem}\label{lem:4.3} If $\vrp: R \to U$ is a supervaluation,

then the set
$$ \langle \vrp(R)\rangle := \vrp(R) \cup e \vrp(R)$$ is a
subsemiring of $U$ (and hence a supertropical semiring itself).
\end{lem}
This can be easily verified.
\begin{defn}\label{defn:4.4} $ $
\begin{enumerate} \ealph
    \item Given supervaluations $\vrp : R \to U$ and $\psi: R \to
    V$ we say that $\vrp$ \textbf{dominates} $\psi$, and write $\vrp \geq
    \psi$, if there exists a semiring homomorphism
$$ \al: \langle \vrp(R)\rangle \to \langle \psi(R)\rangle,$$
necessarily surjective, such that $\psi(a) = \al(\vrp(a))$ for
every $a\in R$. \pSkip

    \item We call $\vrp$ and $\psi$ \textbf{equivalent}, and write
    $\vrp \sim \psi$, of both $\vrp \geq \psi$ and $\psi \geq
    \vrp$.\pSkip

    \item We denote the equivalence class of a supervaluation
    $\vrp$ covering $v$ by $[\vrp]$, and denote the set of all
    these classes by $\Cov(v).$
\end{enumerate}
\end{defn}

We obtain on the set $\Cov(v)$ a partial ordering by declaring
that
$$ [\vrp] \geq [\psi] \quad \text{iff} \quad \vrp \geq \psi.$$
We now have a fairly remarkable fact:

\begin{thm}\label{thm:4.5} The partially ordered set $\Cov(v)$ is a
complete lattice.
\end{thm}

As every complete lattice $\Cov(v)$ has a top element and a bottom
element. The top element can be described explicitly, cf.
\cite[Example 4.5]{IKR1}. The bottom element is the class $[v]$ of
the supervaluation $v: R \to M$, viewed as a supervaluation.
\{N.B. We regard $M$ is as supertropical semifield without
tangible elements.\}

Starting from now, until the end of the paper, \emph{we assume
that  $v: R \to M$ is a strong valuation} (e.g. $R$ is a ring),
and we focus on a particularly good natured class of
supervaluations covering $v$, to be defined as follows.

\begin{defn}\label{defn:4.6}
A supervaluation  $\vrp: R \to U$ covering $v$ is \textbf{strong}
if
$$ \forall a,b \in R: \qquad \vrp(a) + \vrp(b) \in \tT(U) \ \Rightarrow \
\vrp(a + b) = \vrp(a) + \vrp(b).$$
\end{defn}

The strong supervaluations turn out to be ``nearly" semiring
homomorphisms in the GS-sense. More precisely

\begin{prop}\label{prop:4.7}
A supervaluation $\vrp: R \to U$ covering $v$ is strong iff for
all $a, b \in R$
$$ \vrp(a) + \vrp(b) \lmod \vrp(a + b).$$
\end{prop}

We call a supervaluation $\vrp: R \to U$ \textbf{tangible} if all
its values are tangible or zero; $\vrp(R) \subset \tT(U) \cup \{ 0
\}$.

In the next section the strong valuations which are also tangible
will play a useful role. We quote the following important fact, to
be found in \cite[\S11]{IKR1}.
\begin{thm}\label{thm:4.8}
The subset $\Cov_{t,s} (v)$ of $\Cov(v)$ consisting of all classes
$[\vrp] \in \Cov(v)$ with $\vrp$ tangible and strong is a complete
sublattice of $\Cov(v)$. In particular $\Cov_{t,s} (v)$ is not
empty.
\end{thm}

Again the top and the bottom elements of $\Cov_{t,s} (v)$ can be
described explicitly, cf. \cite[Theorem 11.8 and Example
10.16]{IKR1}.

\section{A supertropical version of Kapranov's lemma}

Assume that $R$ is a semiring, $\vrp: R \to U$ is  a strong
supervaluation covering $v: R \to M$, and $\lm = (\lm_1, \dots,
\lm_n)$ a set of variables. We start  out to extend $\vrp$ to a
supervaluation on the polynomial semiring $R[\lm]$ in various
ways.

We first extend $\vrp$ to a map
$$ \tvrp: R[\lm] \to U[\lm]$$ by the formula
$$ \tvrp\bigg(\sum_i c_i \lm ^i \bigg) : = \sum_i \vrp(c_i)\lm^i.$$ Here  we
use that standard monomial notation: $i$ runs through the set of
tuples $i = (i_1,\dots, i_n)$ with  $i_1,\dots, i_n$ in $\bbN_0$;
$\lm^i$ means $\lm_i^{i_1}\cdots \lm_n^{i_n}$; only finitely  many
$c_i$ are not zero. In the same  way we have a map $\tv: R[\lm]
\to M[\lm]$,
$$ \tv\bigg(\sum_i c_i \lm ^i \bigg) : = \sum_i v(c_i)\lm^i.$$
Now we choose a tuple $a = (a_1,\dots,a_n)$ in $R^n$.

It gives us tuples $\vrp(a) = (\vrp(a_1), \dots, \vrp(a_n) )$ in
$U^n$ and $v(a) = (v(a_1), \dots, v(a_n) )$ in $M^n$. Associated
with these tuples we obtain evaluation maps
$$ \ep_a : R[\lm] \to R, \qquad \ep_{\vrp(a)}: U[\lm] \to U, \qquad
\ep_{v(a)}: M[\lm] \to M,$$ by inserting the tuples for the
variables into the polynomials. For example
$$ \ep_{\vrp(a)}\bigg(\sum_i \gm_i \lm ^i \bigg) : = \sum_i \gm_i \vrp(a)^i \qquad (\gm_i \in U).$$
These maps are semiring homomorphisms.

It is then fairly  obvious that the map
$ v \circ \ep_a: R[\lm] \to M$
 is a valuation and
$ \vrp \circ \ep_a: R[\lm] \to U$
 is a supervaluation covering $v
\circ \ep_a$. With some work it can be seen that also
$ \ep_{v(a)} \circ \tv: R[\lm] \to M$
 is a valuation and
$ \ep_{\vrp(a)} \circ \tv: R[\lm] \to U$
 is a supervaluation
covering $\ep_{v(a)} \circ \tv$. \{Here it is important to assume
that $\vrp$ is strong.\}\footnote{The valuations and
supervaluations on $R[\lm]$ ocuring here are again strong, but
this will not matter for the following.}

Now most often the diagram
$$\xymatrix{
           R[\lm] %
            \ar[d]_{\tvrp}     \ar[rr]^{\ep_a}    &&      R
       \ar[d]^{\vrp}
       \\
      U[\lm]    \ar[rr]^{\ep_{\vrp(a)}}   &&        U
 }$$
does not commute. Instead we have

\begin{thm}\label{thm:5.1} \cite[\S13]{IKR1} For any $f \in R[\lm]$
$$ \ep_{\vrp(a)}  \tvrp(f) \lmod \vrp \ep_a(f).$$
\end{thm}
The theorem says in more imaginative terms that the supervaluation
$ \ep_{\vrp(a)}  \tvrp$ is a perturbation of $\vrp \ep_a$ by
noise.

Theorem \ref{thm:5.1} has close relation to an initial key
observation of tropical geometry, Kapranov's Lemma. Let us briefly
indicate what is says.

Assume that $R$ is a field and $f= \sum c_i \lm^i$ is a polynomial
over $R$. It gives us the hypersurface
$$Z(f) := \{ a \in R^n \ | \ f(a) = 0\}.$$
In tropical geometry one relates $Z(f)$ to the so called ``corner
locus", or ``tropical hypersurface", of the polynomial
$$ \tv(f) = \sum_i v(c_i) \lm ^i \in M[\lm].$$
Notice that, if a tuple $\xi \in M^n$ is given, then
$$ \tv(f)(\xi) = \max_i( v(c_i)\xi^i).$$
The corner locus $Z_0(\tv(f))$ is defined as the set of all tuples
$\xi \in M^n$, where this maximum is attained at least at two
indices. Kapranov's lemma states that
$$ v(Z(f)) \subset Z_0(\tv(f)),$$
cf. \cite[Lemma 2.1.4]{EKL}.

It can be deduced from Theorem \ref{thm:5.1} as follows.

We choose a tangible strong supervaluation $\vrp: R \to U$
covering $v$, which is possible by Theorem \ref{thm:4.8}. Let $a
\in Z(f)$. Then $\vrp \ep_a(f) = \vrp(f(a)) =0$. Theorem
\ref{thm:5.1} tells us that
$$ \ep_{\vrp(a)}  \tvrp(f) =
\sum_i \vrp(c_i) \vrp(a)^i\lmod 0,$$ i.e., this sum is ghost. But
each summand $\vrp(c_i) \vrp(a)^i$ is tangible or zero. From the
law ST3 in \S3 we infer that the maximum of the values
$$e \vrp(c_i) \vrp(a)^i = (e\vrp(c_i))(e \vrp(a)^i) = v(c_i) v(a)^i$$
is attained by more than one index. In other words, $v(a)$ is an
element of the corner locus $Z_0(\tv(f))$. Thus indeed $ v(Z(f))
\subset Z_0(\tv(f)).$

Theorem \ref{thm:5.1} says more than the classical Kapranov lemma,
not only since a semiring $R$ instead of a field $R$ is admitted,
but also since it contains a statement about points $a \in R^n$
with $f(a) \neq 0$.

Finally, if $\vrp$ and $\psi$ are strong supervaluations covering
$v$ with $\vrp \geq \psi$, the statement of Theorem \ref{thm:5.1}
for $\vrp$ formally implies the same statement for $\psi$. Thus
Theorem~\ref{thm:5.1} seems to be ``best", if $\vrp$ is the top
element of $\Cov_{t,s}(v)$, at least if we focus on tangible
supervaluations.

To exploit all this, more work will be needed than what has been
done in \cite{IKR1}.

\end{document}